\documentclass{article}
\usepackage{amsmath,amsthm,amscd,amssymb} 

 \newtheorem{theorem}{Theorem}[section]

 \newtheorem{proposition}[theorem]{Proposition}
 \newtheorem{lemma}[theorem]{Lemma}
 \newtheorem{corollary}[theorem]{Corollary}
 \newtheorem{remark}[theorem]{Remark}
 \newtheorem{definition}[theorem]{Definition}
 
 \newtheorem{example}[theorem]{Example}

 \numberwithin{equation}{section}

 \def\subrel#1#2{\mathrel{\mathop{#2}\limits_{#1}}}
 \def\qqq{\qed\vspace{2mm}} 

 \def\fp#1{{\mathbb F}_{\! p^{#1}}}

 \def\z2{{\mathbb Z}/2}

 \def\hz2{H\z2}

 \def\power#1{[\![#1]\!]}

 \def\gal#1{\mbox{\rm Gal}(#1)}

 \def\inverselimit#1{\ \subrel{\subrel{#1}{\longleftarrow}}{\lim}}
 \def\directlimit#1{\ \subrel{\subrel{#1}{\longrightarrow}}{\lim}}

\begin{document}

\title{Equivariance of generalized Chern characters}
\author{Takeshi Torii}
\date{\footnote{
Date: March 27, 2009.}
\setcounter{footnote}{-1}
\footnote{
\textup{2000} {\it Mathematics Subject Classification}.
Primary 55N22 ; Secondary 55P43, 55S05.}
\setcounter{footnote}{-1}
}

\maketitle

\vspace{-15mm}

\begin{abstract}

In this note 
some generalization 
of the Chern character is discussed
from the chromatic point of view.
We construct a multiplicative $G_{n+1}$-equivariant 
natural transformation $\Theta$ 
from some height $n+1$ cohomology theory $E^*(-)$
to the height $n$ cohomology theory 
$K^*(-)\widehat{\otimes}_{\mathbf F}L$,
where $K^*(-)$ is essentially the $n$th Morava $K$-theory.  
As a corollary, it is shown that the $G_n$-module $K^*(X)$ 
can be recovered from the $G_{n+1}$-module $E^*(X)$.
We also construct a lift of $\Theta$ to a natural transformation 
between characteristic zero cohomology theories.  

\end{abstract}

\section{Introduction}

%

In the stable homotopy category ${\mathcal S}$ of $p$-local spectra,
there is a filtration of full subcategories ${\mathcal S}_n$,
where the objects of ${\mathcal S}_n$ consist of $E(n)$-local spectra.
The difference of the each step of this filtration
is equivalent to the $K(n)$-local category.
So it can be considered that the stable homotopy category ${\mathcal S}$
is built up from $K(n)$-local category.
In fact, the chromatic convergence theorem (cf. \cite{Ravenel-2})
says that the tower
$\cdots\to L_{n+1}X\to L_n X\to\cdots\to L_0X$ 
recovers a finite spectrum $X$, that is,  
$X$ is homotopy equivalent to the
homotopy inverse limit of the tower.
Furthermore,
the chromatic splitting conjecture (cf. \cite{Hovey}) implies that
the $p$-completion of a finite spectrum $X$ is a direct summand of 
the product $\prod_n L_{K(n)}X$.
This means that 
it is not necessarily to reconstruct the tower but
it is sufficient to know all $L_{K(n)}X$
to obtain some information of $X$.

The weak form of the chromatic splitting conjecture means that
the canonical map 
$L_n (S^0)^{\wedge}_p \to L_n L_{K(n+1)} S^0$
is a split monomorphism,
where $S^0$ is the sphere spectrum and
$(S^0)^{\wedge}_p$ is its $p$-completion.  
In \cite[Remark~3.1.(i)]{Minami-1}
Minami indicated that  
the weak form of the chromatic splitting conjecture 
implies that there is a natural map $\rho$ for a finite spectrum $X$
from the $K(n+1)$-localization $L_{K(n+1)}X$ 
to the $K(n)$-localization $L_{K(n)}X$ such that
the following diagram commutes:
\begin{equation}\label{topological-diagram}
\begin{array}{rcl}
      &  X    &  \\
      & {\begin{minipage}{1.2cm}
        \begin{picture}(30,20) 
           \put(-38,15){${\scriptstyle{\eta_{K(n+1)}}}$}
           \put(28,23){\vector(1,-1){20}}
           \put(43,15){${\scriptstyle{\eta_{K(n)}}}$}
           \put(5,23){\vector(-1,-1){20}} 
        \end{picture}
        \end{minipage}} &  \\
    L_{K(n+1)}X& 
    \stackrel{\rho}{\hbox to 2cm{\rightarrowfill}} & 
    L_{K(n)}X,
       \end{array}   \end{equation}  
where $\eta_{K(n)}$ and $\eta_{K(n+1)}$ are the localization maps. 
In this note we would like to consider an algebraic analogue 
of this diagram.

Let ${\mathbf F}$ be an algebraic extension of the 
prime field $\fp{}$ which contains
$\fp{n}$ and $\fp{n+1}$.
Let $E_n$ be the Morava $E$-theory
with the coefficient ring 
$W\power{u_1,\ldots,u_{n-1}}[u^{\pm 1}]$
where $W=W({\mathbf F})$ is the ring of Witt vectors
with coefficients in ${\mathbf F}$.
The group of multiplicative cohomology operations 
on $E_n$ is the extended Morava stabilizer group 
$G_n=\Gamma\ltimes S_n$,
where $\Gamma=\gal{{\mathbf F}/\fp{}}$ and
$S_n$ is the $n$th Morava stabilizer group.
Let ${\mathcal M}_n$ be the category of 
twisted $E_{n*}$-$G_n$-modules.
Then $E_n$-(co)homology gives a functor from 
${\mathcal S}$ to ${\mathcal M}_n$, and 
$E_n$-(co)homology is considered to be an algebraic 
approximation of the localization map
$\eta_{K(n)}: X\to L_{K(n)}X$.
So an algebraic analogue of the problem to construct
the diagram (\ref{topological-diagram}) is 
the following:
Is there an algebraic functor 
$\mu: {\mathcal M}_{n+1}\to {\mathcal M}_n$
such that the following diagram commutes?
\begin{equation}\label{algebraic-diagram}
\begin{array}{rcl}
      & {\mathcal S}    &  \\
      & {\begin{minipage}{1.2cm}
        \begin{picture}(30,20) 
           \put(-27,15){${\scriptstyle{E_{n+1}}}$}
           \put(28,23){\vector(1,-1){20}}
           \put(43,15){${\scriptstyle{E_n}}$}
           \put(5,23){\vector(-1,-1){20}} 
        \end{picture}
        \end{minipage}} &  \\
    {\mathcal M}_{n+1} & 
    \stackrel{\mu}{\hbox to 2cm{\rightarrowfill}} & 
    {\mathcal M}_n.
       \end{array}   
\end{equation}  


We recall the classical Chern characters.
The Chern character is a multiplicative natural transformation
from $K$-theory to the rational cohomology:
\[ ch: K\longrightarrow H{\mathbf Q}[w^{\pm 1}]\cong
       \prod_{i\in {\mathbf Z}}\Sigma^{2i}H{\mathbf Q}. \]
The formal group law associated with $K$-theory
is the multiplicative formal group law.
So its height is one.
On the other hand,
the height of the formal group law associated with the rational cohomology
is considered to be zero.
So the Chern character is considered to be a transformation
from a height $1$ theory to a height zero theory.
The map $\mu$ in (\ref{algebraic-diagram})
is a transformation from the height $n+1$
theory $E_{n+1}$ to the height $n$ theory $E_n$.
So $\mu$ should be a generalization of Chern character
in some sense.
Such a generalized Chern character
have been constructed and studied by Ando, Morava and Sadofsky \cite{AMS}.

There is a modulo $I_n$-version of the problem
to construct the diagram (\ref{algebraic-diagram}),
where $I_n$ is the invariant prime ideal $(p,v_1,\ldots,v_{n-1})$.
Let $E^*(-)$ and $K^*(-)$ be the complex oriented cohomology theories
with coefficient rings $E_*={\mathbf F}\power{u_n}[u^{\pm 1}]$
and $K^*(-)={\mathbf F}[w^{\pm 1}]$, respectively,
where $u^{-(p^{n+1}-1)}=v_{n+1}$ and 
$u_n u^{-(p^n-1)}=v_n= w^{-(p^n-1)}$.
We denote by ${\mathcal M}'_{n+1}$ (resp. ${\mathcal M}'_n$)
the category of twisted $E_*$-$G_{n+1}$-modules
(resp. $K_*$-$G_n$-modules).
Then the modulo $I_n$-version of the problem
to construct the diagram (\ref{algebraic-diagram})
is as follows:
Is there an algebraic functor 
$\mu': {\mathcal M}'_{n+1}\to {\mathcal M}'_n$
such that the following diagram commutes?
\begin{equation}\label{algebraic-diagram-mod-In}
\begin{array}{rcl}
      & {\mathcal S}    &  \\
      & {\begin{minipage}{1.2cm}
        \begin{picture}(30,20) 
           \put(-17,13){${\scriptstyle{E}}$}
           \put(28,23){\vector(1,-1){20}}
           \put(43,13){${\scriptstyle{K}}$}
           \put(5,23){\vector(-1,-1){20}} 
        \end{picture}
        \end{minipage}} &  \\
    {\mathcal M}'_{n+1} & 
    \stackrel{\mu'}{\hbox to 2cm{\rightarrowfill}} & 
    {\mathcal M}'_n.
       \end{array}   
\end{equation}  


In \cite{Torii} we have studied the relationship between 
the formal group laws $F_{n+1}$ and $H_n$ associated
with $E^*(-)$ and $K^*(-)$, respectively.
There is 
a totally ramified Galois extension $L$ of infinite degree
over the fraction field ${\mathbf F}((u_n))$ of $E_0$,
and there is an isomorphism between $F_{n+1}$ and $H_n$ over $L$.
We have shown that 
the pro-finite group ${\mathcal G}=\Gamma\ltimes (S_n\times S_{n+1})$
acts on $(F_{n+1},L)\cong (H_n,L)$.
 
The following is the main theorem of this note.

\begin{theorem}[Theorem~\ref{main-theorem}]\label{intro-main}
Let $p$ be an odd prime.
Then there is a $G_{n+1}$-equivariant multiplicative 
stable cohomology operation
\[ \Theta: E^*(-)\longrightarrow K^*(-)\widehat{\otimes}_{\mathbf F}L \]
such that $\Theta$ induces a natural isomorphism of 
${\mathcal G}=\Gamma\ltimes (S_n\times S_{n+1})$-modules:
\[ E^*(X)\widehat{\otimes}_{E_0}L\stackrel{\cong}{\longrightarrow}
   K^*(X)\widehat{\otimes}_{\mathbf F}L \]
for all spectra $X$.
\end{theorem} 

On the right hand side of the isomorphism in Theorem~\ref{intro-main},
the subgroup $S_{n+1}$ of ${\mathcal G}$ acts on $L$ only and
its invariant ring is the subfield ${\mathbf F}$:
$H^0(S_{n+1};L)={\mathbf F}$.
This implies the following corollary.

\begin{corollary}[Corollary~\ref{main-corollary}]
There are natural isomorphisms of $G_n$-modules:
\[ \begin{array}{rcl}
    K^*(X)&\cong& H^0(S_{n+1};E^*(X)\widehat{\otimes}_{E}
    L),\\[2mm]
    K_*(X)&\cong& H^0(S_{n+1};E_*(X){\otimes}_{E}
    L),\\
   \end{array} \]
for all spectra $X$.
If $X$ is a space, then 
these are also isomorphisms of graded commutative rings.
\end{corollary}

This corollary gives us an answer of the 
modulo $I_{n-1}$-version of the problem.
We define 
$\mu'(M)=H^0(S_n;M\otimes_{E_*}L[u^{\pm 1}])$ for 
a twisted $E_*$-$G_{n+1}$-module $M$.
Then we obtain a functor 
from the category of twisted $E_*$-$G_{n+1}$-modules
to the category of twisted $K_*$-$G_n$-modules:
$\mu': {\mathcal M}'_n  \to {\mathcal M}'_{n-1}$,
which makes the triangle (\ref{algebraic-diagram-mod-In}) 
commutative.

Furthermore,
we can lift $\Theta$
to a natural transformation between 
characteristic $0$ cohomology theories.
There is a complete 
discrete valuation ring $T$ of characteristic $0$
with uniformizer $p$ and residue fields $L$.
We regard $T\power{w_i}=T\power{w_1,\ldots,w_{n-1}}$
as an $E_n=W\power{w_1,\ldots,w_{n-1}}$-algebra by
obvious way.
Also, we can regard $T\power{u_i}=T\power{u_1,\ldots,u_{n-1}}$
as an $E_{n+1}=W\power{u_1,\ldots,u_n}$-algebra.

\begin{theorem}[Theorem~\ref{Characteristic-zero}]
There is a $G_{n+1}$-equivariant multiplicative stable operation
\[ ch: E_{n+1}^*(-)\longrightarrow E_n^*(-)
       \widehat{\otimes}_{E_n} T\power{w_i}, \]
such that this induces a natural isomorphism of ${\mathcal G}$-modules:
\[ E_{n+1}^*(X)\widehat{\otimes}_{E_{n+1}} T\power{u_i}
   \stackrel{\cong}{\longrightarrow}
   E_n^*(X)\widehat{\otimes}_{E_n} T\power{w_i} \]
for all spectra $X$.
\end{theorem}

The organization of this note is as follows:
In \S\ref{FGL}
we review the Lubin-Tate's deformation theory of formal group laws
and the results of \cite{Torii} on the degeneration of formal group laws.
In \S\ref{stable-operations}
we study the relationship between the stable natural transformations of 
even-periodic complex oriented cohomology theories
and the homomorphisms of their formal group laws.
In \S\ref{generalized-Chern-ch}
we construct a multiplicative $G_{n+1}$-equivariant 
natural transformation $\Theta$ from $E^*(-)$ to 
$K^*(-)\widehat{\otimes}L$
and prove the main theorem.
In \S\ref{lift-to-ch-0}  
we construct a lift of $\Theta$ to a natural transformation
of characteristic zero cohomology theories.

\section{Formal group laws}\label{FGL}

In this section we review the deformation theory of
formal group laws.
In the following of this note a formal group law means 
a one-dimensional commutative formal group law.


Let $R_1$ and $R_2$ be two (topological) commutative rings.
Let $F_1$ (resp. $F_2$) be a 
formal group law over $R_1$ (resp. $R_2$).
We understand that a homomorphism from $(F_1,R_1)$ 
to $(F_2,R_2)$ 
is a pair $(f,\alpha)$ of a (continuous) ring homomorphism
$\alpha:R_2\to R_1$ and a homomorphism 
$f:F_1\to \alpha^*F_2$ in the usual sense,
where $\alpha^*F_2$ is the formal group law obtained from $F_2$
by the base change induced by $\alpha$.
We denote the set of all such pairs by 
\[ \mathbf{FGL}((F_1,R_1), (F_2,R_2)) .\]
If $R_1$ and $R_2$ are topological rings,
then we denote the subset of $\mathbf{FGL}((F_1,R_1),(F_2,R_2))$
consisting of $(f,\alpha)$ such that $\alpha$ is continuous by
\[ \mathbf{FGL}^c((F_1,R_1), (F_2,R_2)) .\]
The composition of two homomorphisms $(f,\alpha):(F_1,R_1)\to (F_2,R_2)$
and $(\beta,g):(F_2,R_2)\to (F_3,R_3)$ is defined as 
$(\alpha^*g\circ f,\alpha\circ\beta):(F_1,R_1)\to (F_3,R_3)$:
\[ F_1\stackrel{f}{\longrightarrow}\alpha^*F_2
      \stackrel{\alpha^*g}{\longrightarrow}
      \alpha^*(\beta^*F_3)=(\alpha\circ\beta)^*F_3.\]
A homomorphism $(f,\alpha):(F_1,R_1)\to (F_2,R_2)$ is an isomorphism
if there exists a homomorphism $(g,\beta):(F_2,R_2)\to (F_1,R_1)$ 
such that $(f,\alpha)\circ (g,\beta)=(X, id)$ and
$(g,\beta)\circ (f,\alpha)=(X, id)$.
Then a homomorphism $(f,\alpha):(F_1,R_1)\to (F_2,R_2)$ is an
isomorphism if and only if $\alpha$ is a (topological) ring isomorphism
and $f$ is an isomorphism in the usual sense. 

There is a $p$-typical formal group law $H_n$
over the prime field $\fp{}$ with $p$-series
\[ [p]^{H_n}(X)=X^{p^n},\]
which is called the height $n$ Honda
formal group law.
Let ${\mathbf F}$ be an algebraic extension of 
the finite field $\fp{n}$ with $p^n$ elements, 
and we suppose that  $H_n$ is defined over ${\mathbf F}$.
The automorphism group $S_n$ of $H_n$ over 
${\mathbf F}$ in the usual sense
is the $n$th Morava stabilizer group $S_n$, 
which is isomorphic to the unit group of
the maximal order of the central division algebra
over the $p$-adic number field ${\mathbf Q}_p$
with invariant $1/n$.
We denote by $G_n$ the automorphism group of $H_n$
over ${\mathbf F}$ in the above sense:
\[ G_n=\mathbf{Aut}(H_n, {\mathbf F}).\]
Then the following lemma is well-known.

\begin{lemma}
The automorphism group $G_n$ is isomorphic to
the semi-direct product $\Gamma\ltimes S_n$, 
where $\Gamma$ is the Galois group $\gal{{\mathbf F}/\fp{}}$.
\end{lemma}

We recall Lubin and Tate's
deformation theory of formal group laws \cite{Lubin-Tate}.
Let $R$ be a complete Noetherian local ring with 
maximal ideal $I$ such that
the residue field $k=R/I$ is of characteristic $p>0$.
Let $G$ be a formal group law over $k$ of height $n<\infty$.
Let $A$ be a complete Noetherian local $R$-algebra with maximal
ideal ${\mathfrak m}$. 
We denote by $\iota$ the canonical inclusion of residue fields
$k\subset A/{\mathfrak m}$ induced by the $R$-algebra structure.
A deformation of $G$ to $A$ is a 
formal group law $\widetilde{G}$ over $A$
such that $\iota^*G=\pi^*\widetilde{G}$
where $\pi:A\to A/{\mathfrak m}$ is the canonical projection.
Let $\widetilde{G}_1$ and $\widetilde{G}_2$ 
be two deformations of $G$ to $A$.
We define a $*$-isomorphism between $\widetilde{G}_1$ and 
$\widetilde{G}_2$ as an isomorphism
$\widetilde{u}:\widetilde{G}_1\to\widetilde{G}_2$ over $A$
such that $\pi^*\widetilde{u}$ is the identity map
between $\pi^*\widetilde{G}_1=\iota^*G=\pi^*\widetilde{G}_2$.
Then it is known that there is at most one $*$-isomorphism
between $\widetilde{G}_1$ and $\widetilde{G}_2$.
We denote by ${\mathcal C}(R)$ the category of complete Noetherian
local $R$-algebras with local $R$-algebra homomorphisms as morphisms. 
For an object $A$ of ${\mathcal C}(R)$,
we let $\mathbf{DEF}(A)$ be the set of all $*$-isomorphism classes of 
the deformations of $G$ to $A$.
Then $\mathbf{DEF}$ defines a functor from ${\mathcal C}(R)$
to the category of sets.
Let $R\power{t_i}=R\power{t_1,\ldots,t_{n-1}}$ be a formal power series ring 
over $R$ with $n-1$ indeterminates.
Note that $R\power{t_i}$ is an object of ${\mathcal C}(R)$.
Lubin and Tate constructed a formal group law 
$F(t_i)=F(t_1,\ldots,t_{n-1})$ over $R\power{t_i}$
such that 
for every deformation $\widetilde{G}$ of $G$ to $A$,
there is a unique local $R$-algebra homomorphism 
$\alpha:R\power{t_i}\to A$
such that $\alpha^*F(t_i)$
is $*$-isomorphic to $\widetilde{G}$. 
Hence 
the functor $\mathbf{DEF}$ is represented by
$R\power{t_i}$:
\[ \mathbf{DEF}(A)\cong \mbox{\rm Hom}_{{\mathcal C}(R)}
   (R\power{t_i},A)\]
and $F(t_i)$ is a universal object.


\begin{lemma}\label{key-lemma-ch0}
Let $F$ and $G$ be formal group laws of height $n<\infty$
over a field $k$ of characteristic $p>0$ and  
$(\overline{f},\overline{\alpha})$ 
an isomorphism from $(F,k)$ to $(G,k)$.
Let $R$ be a complete Noetherian local ring with residue field $k$
and $\alpha$ a ring automorphism of $R$ such that
$\alpha$ induces $\overline{\alpha}$ on the residue field.
Let $\widetilde{F}$ (resp. $\widetilde{G}$)
be a universal deformation of $F$ (resp. $G$)
over $R\power{u_i}=R\power{u_1,\ldots,u_{n-1}}$
(resp. $R\power{w_i}=R\power{w_1,\ldots,w_{n-1}}$).
Then there is a unique isomorphism $(g,\beta)$
from $(\widetilde{F},R\power{u_i})$ to $(\widetilde{G},R\power{w_i})$
such that 
$(g,\beta)$ induces $(\overline{f},\overline{\alpha})$
on the residue field
and $i\circ \alpha=\beta \circ j$,
where $i:R\to R\power{u_i}$ and $j:R\to R\power{w_i}$
are canonical inclusions. 
\end{lemma}

\proof
First, we show that there is such a homomorphism.
Let $f(X)\in R[X]$ be a lift of $\overline{f}(X)\in k[X]$
such that $f(0)=0$. 
Set $F'(X,Y)=f(\widetilde{F}(f^{-1}(X),f^{-1}(Y)))$.
Then $(F',R\power{u_i})$ is a deformation of
$\overline{\alpha}^*G$. 
We denote by $R'\power{u_i}$
the ring $R\power{u_i}$ with the $R$-algebra structure
given by $R\stackrel{\alpha}{\to}R\stackrel{i}{\to}R\power{u_i}$.
Then $(F', R'\power{u_i})$ is a deformation of $G$.
Since $\widetilde{G}$ is a universal deformation of $G$,
there exists a continuous $R$-algebra homomorphism 
$\beta: R\power{w_i}\to R'\power{u_i}$
and a $*$-isomorphism $\widetilde{u}: F'\to \beta^*\widetilde{G}$.
Then $(g,\beta)=(\widetilde{u}\circ f, \beta): 
(\widetilde{F}, R'\power{u}_i)\to (\widetilde{G}, R\power{w_i})$
is a lift of $(\overline{f},\overline{\alpha}): (F,k)\to (G,k)$.

By the same way,
we can construct a lift $(h,\gamma)$ of 
$(\overline{f},\overline{\alpha})^{-1}$.
Then $(h,\gamma)\circ (g,\beta)$ is a lift of 
$(X,id): (F,k)\to (F,k)$.
Note that $\beta\circ\gamma: R\power{u_i}\to R\power{u_i}$
is a continuous $R$-algebra homomorphism.   
Since $(\widetilde{F}, R\power{u_i})$ is a universal deformation,
$(h,\gamma)\circ (g,\beta)=(X,id)$ by the uniqueness.
Similarly,
we obtain that $(g,\beta)\circ (h,\gamma)=(X,id)$.
Hence we see that $(g,\beta)$ is an isomorphism
and a unique lift of 
$(\overline{f},\overline{\alpha}): (F,k)\to (G,k)$.
\qqq

Let ${\mathbf F}$ be an algebraic extension of $\fp{}$
which contains $\fp{n}$ and $\fp{n+1}$.
Let $W=W({\mathbf F})$ be the ring of Witt vectors
with coefficients in ${\mathbf F}$.
We define $E_n$ to be
a formal power series ring over $W$
with $(n-1)$ indeterminates:
\[ E_n=W\power{w_1,\ldots,w_{n-1}}.\]
The ring $E_n$ is a complete Noetherian local ring 
with residue field ${\mathbf F}$.
There is a $p$-typical formal group law $\widetilde{F}_n$
over $E_n$ with the $p$-series:
\begin{equation}\label{p-series-tilde-Fn+1}
 [p]^{\widetilde{F}_n}(X)=pX+_{\widetilde{F}_n}w_1X^p
                            +_{\widetilde{F}_n}w_2X^{p^2}
                            +_{\widetilde{F}_n}\cdots
                            +_{\widetilde{F}_n}w_{n-1}X^{p^{n-1}}
                            +_{\widetilde{F}_n}X^{p^n}.
\end{equation}
The formal group law $\widetilde{F}_n$ is a deformation of $H_n$
to $E_n$.  
The following lemmas are well-known.

\begin{lemma}\label{universal-deformation}
$(\widetilde{F}_n, E_n)$ is a universal deformation 
of $(H_n,{\mathbf F})$.
\end{lemma}

\begin{lemma}\label{auto-gn}
The automorphism group $\mathbf{Aut}^c(\widetilde{F}_n,E_n)$ 
is isomorphic to $G_n$.
\end{lemma}

As in $E_n$, we define $E_{n+1}$ to be 
a formal power series ring over $W$
with $n$ indeterminates:
\[ E_{n+1}=W\power{u_1,\ldots,u_n},\]
and there is a universal deformation
$(\widetilde{F}_{n+1}, E_{n+1})$ of 
the height $(n+1)$ Honda group law $(H_{n+1},{\mathbf F})$. 
Let $E=E_{n+1}/I_n={\mathbf F}\power{u_n}$,
where $I_n=(p,u_1,\ldots,u_{n-1})$.
Let $F_{n+1}=\pi^*\widetilde{F}_{n+1}$,
where $\pi$ is the quotient map $E_{n+1}\to E$.
Then $F_{n+1}$ is a deformation of $H_{n+1}$ to $E$.
The following lemma is easy.

\begin{lemma}\label{automorphism_group_F_n+1}
The automorphism group $\mathbf{Aut}^c(F_{n+1},E)$
is isomorphic to $G_{n+1}$.
\end{lemma}

If we suppose that $F_{n+1}$ is defined 
over the quotient field $M={\mathbf F}((u_n))$ of $E$,
then its height is $n$.
Since the formal group laws over a separably closed field
is classified by their height,
there is an isomorphism $\Phi$ between $F_{n+1}$ and $H_n$
over the separable closure $M^{sep}$ of $M$ 
(cf. \cite{Lazard, Hazewinkel}).
We fix such an isomorphism $\Phi$.
Since $\Phi: F_{n+1}\to H_n$ is a homomorphism between $p$-typical formal
group laws, 
$\Phi$ has a following form:
\[ \Phi(X)=\sum_{i\ge 0}{}^{H_n}\Phi_i X^{p^i}.\]
Let $L$ be the extension field of $M$
obtained by adjoining all the coefficients of the isomorphism $\Phi$.
So $(\Phi,id_L)$ is
an isomorphism from $(F_{n+1},L)$ to $(H_n,L)$:
\[   (\Phi, id_L): (F_{n+1},L)\stackrel{\cong}{\longrightarrow}
                   (H_n,L).\]
Note that $L$ is a totally ramified Galois extension
of infinite degree over $M$ with Galois group isomorphic to $S_n$
\cite{Gross, Torii}.
Set ${\mathcal G}=\Gamma\ltimes (S_{n+1}\times S_n)$.
Then $G_{n+1}=\Gamma\ltimes S_{n+1}$ and $G_n=\Gamma\ltimes S_n$
are subgroups of ${\mathcal G}$.
In \cite{Torii} we have shown the following theorem.

%
%
%

\begin{theorem}[{cf. \cite[\S2.4]{Torii}}]\label{FGL-main-theorem-v3}
The pro-finite group ${\mathcal G}$ acts on 
$(F_{n+1},L)\stackrel{}{\cong} (H_n,L)$.
The action of the subgroup $G_{n+1}$ on $(F_{n+1},L)$
is an extension of the action on $(F_{n+1},E)$,
and the action of the subgroup $G_n$ on $(H_n,L)$
is an extension of the action on $(H_n,{\mathbf F})$.
\end{theorem}

\section{Stable operations of cohomology theories}\label{stable-operations}


In this section we recall and study the stable cohomology operations
between Landweber exact cohomology theories over $P(n)$.
The treatment is standard as in \cite{Adams, Miller-2, Wurgler-2}.   

For a spectrum $h$,
we denote by
$h^*(-)$ (resp. $h_*(-)$) the associated 
generalized cohomology (resp. homology) theory.
For spectra $h$ and $k$,
we denote by ${\mathbf C}(h,k)$ 
(resp. ${\mathbf H}(h,k)$)
the set of all degree $0$
stable cohomology (resp. homology) operations 
from $h$ to $k$.
Then ${\mathbf C}(h,k)$ is naturally identified with
the set of all degree $0$ morphisms from $h$ to $k$
in the stable homotopy category.
There is a natural surjection from 
${\mathbf C}(h,k)$ to ${\mathbf H}(h,k)$
and the kernel consists of phantom maps
(cf. \cite[Chapter~4.3]{Margolis}).

We say that a graded commutative ring $h_*$ is even-periodic 
if there is a unit $u\in h_2$ of degree $2$
and $h_{odd}=0$.
Note that $h_*=h_0[u^{\pm 1}]$ if $h_*$ is even-periodic.
We say that a ring spectrum $h$ is even-periodic
if the coefficient ring $h_*$ is even-periodic.

\begin{definition}\rm
Let $R$ be a 
commutative ring.
A topological $R$-module $M$ is said to be 
{\it linearly topologized} if $M$ has a fundamental 
neighbourhood system at the zero consisting 
of the open submodules.
A linearly topologized $R$-module $M$ 
is said to be {\it linearly compact}
if it is Hausdorff and it has the finite intersection
property with respect to the closed cosets
A topological ring $R$ is linearly compact
if $R$ is linearly compact as an $R$-module.
(cf. 
\cite[Definition~2.3.13]{Hovey-Palmieri-Strickland}).
\end{definition}

\begin{example}\rm
A linearly topologized compact Hausdorff (e.g. profinite) 
module is linearly compact.
If $R$ is a complete Noetherian local ring,
then a finitely generated $R$-module is linearly compact.
In particular,
a finite dimensional vector space over a field 
is linearly compact. 
\end{example}

\begin{lemma}[cf. {\cite[Corollary~2.3.15]{Hovey-Palmieri-Strickland}}]
\label{exact-inverse-limit}
Let ${\mathcal I}$ be a filtered category.
The inverse limit functor indexed by ${\mathcal I}$
is exact
in the category of linearly compact modules
and continuous homomorphisms.
\end{lemma}


For a spectrum $X$,
we denote by $\Lambda(X)$ the category
whose objects are maps $Z\stackrel{u}{\to}X$ 
such that $Z$ is finite,
and whose morphisms are maps $Z\stackrel{v}{\to}Z'$
such that $u'v=u$.
Then $\Lambda(X)$ is an essentially small filtered category.

\begin{lemma}\label{no-phantom}
If $k$ is even-periodic, and
$k_0$ is Noetherian and linearly compact,
then there is no phantom maps to $k$.
\end{lemma}

\proof
For a finite spectrum $Z$,
$k^0(Z)$ is a finitely generated module over $k^0$,
and hence $k^0(Z)$ is linearly compact.
By Lemma~\ref{exact-inverse-limit},
$k^0(X)\cong\ \subrel{\longleftarrow}{\lim}k^0(Z)$,
where the inverse limit is taken over $\Lambda(X)$.
This means that there is no phantom maps to $k$.
\qqq

\begin{corollary}\label{homology-operation}
Suppose that a spectrum $k$ is even-periodic, and
$k_0$ is Noetherian and linearly compact.
Then the natural map 
${\mathbf C}(h,k)\to {\mathbf H}(h,k)$ is an isomorphism.  
\end{corollary}       

\proof
Since ${\mathbf C}(h,k)\to {\mathbf H}(h,k)$
is surjective and the kernel consists of phantom maps,
the corollary follows from Lemma~\ref{no-phantom}.
\qqq

\begin{definition}\rm
We denote by $\mathbf{Mult}(h,k)$ the set of
all multiplicative stable cohomology operations 
from $h^*(-)$ to $k^*(-)$.
If $h^*(-)$ and $k^*(-)$ have their values in the category 
of linear compact modules,
then we denote by $\mathbf{Mult}^c(h,k)$ 
the subset of $\mathbf{Mult}(h,k)$ consisting
of $\theta$ such that $\theta: h^*(X)\to k^*(X)$ 
is continuous for all $X$.
\end{definition}

If $h^*(-)$ is a complex oriented cohomology theory,
then the orientation class $X_h\in h^2({\mathbf C}P^{\infty})$
gives a formal group law $F$ of degree $-2$.
Furthermore, if $h$ is even-periodic,
then a unit $u\in h_2$ gives a degree $0$ formal group law
by $F_h(X,Y)= u F(u^{-1}X, u^{-1}Y)$.
In the following of this section
we suppose that $h^*(-)$ and $k^*(-)$ are complex orientable
and even-periodic.
Furthermore, we fix a unit $u\in h_2$ (resp. $v\in k_2$)
and an orientation class $X_h\in h^2({\mathbf C}P^{\infty})$
(resp. $X_k\in k^2({\mathbf C}P^{\infty})$). 
Then we obtain a degree $0$ formal group law
$F_h$ (resp. $F_k$) associated with $h$ (resp. $k$)
as above.
A multiplicative cohomology operation $\theta: h^*(-)\to k^*(-)$
gives a ring homomorphism $\alpha: h_0\to k_0$ and
an isomorphism $f: F_k\to \alpha_*F_h$ of 
formal group laws.
Note that $f(X)=\theta(u)\widetilde{f}(v^{-1}X)$,
where $\widetilde{f}(X_k)=\theta(X_h)$.
In particular,
$f'(0)=\theta (u) v^{-1}$ is a unit of $k_0$.
Hence we obtain a map from $\mathbf{Mult}(h,k)$
to $\mathbf{FGL}((F_k,k_0),(F_h,h_0))$:
\[ \Xi:  \mathbf{Mult}(h,k) \to \mathbf{FGL}((F_k.k_0),(F_h,h_0)) .\] 
If $h_0$ and $k_0$ are Noetherian and linearly compact,
then $h^*(X)$ and $k^*(X)$ are linearly compact modules.
Then we see that $\Xi$ induces 
\[ \Xi^c: \mathbf{Mult}^c(h,k) \to \mathbf{FGL}^c((F_k.k_0),(F_h,h_0)) .\]   

\begin{remark}\label{cotinuity-remark}\rm
Let $h_0$ and $k_0$ be Noetherian and linearly compact. 
If $\theta\in\mathbf{Mult}(h,k)$ induces a continuous
ring homomorphism $h_0\to k_0$,
then $\theta\in\mathbf{Mult}^c(h,k)$.
\end{remark}


Let $p$ be a prime number and 
$BP$ the Brown-Peterson spectrum at $p$.
There is a $BP$-module spectrum $P(n)$
with coefficient $P(n)_*=\fp{}[v_n,v_{n+1},\ldots]$. 
If $p$ is odd,
then $P(n)$ is a commutative $BP$-algebra spectrum.
As usual, we set $P(0)=BP$.

Let $R_*$ be a graded commutative ring over ${\mathbf Z}_{(p)}$.
We suppose that there is a $p$-typical 
formal group law $F$ of degree $-2$ over $R_*$.
Since the associated formal group law to $BP$ 
is universal with respect to $p$-typical ones,
there is a unique ring homomorphism $r: BP_*\to R_*$.
We suppose that $r(v_i)=0$ for $0\le i<n$.
Then we obtain a ring homomorphism $\widetilde{r}: P(n)_*\to R_*$.
The functor $R_*\otimes_{P(n)_*} P(n)_*(-)$ is 
a generalized homology theory 
and $R_*\otimes_{P(n)_*}P(n)^*(-)$
is a cohomology theory on the category of finite spectra
if $v_n,v_{n+1},\ldots$ is a regular sequence in $R_*$
by the exact functor theorem~\cite{Landweber, Yagita}.
We say that such a graded ring $R_*$ is
Landweber exact over $P(n)_*$.
For a spectrum $X$, we define 
\[ R_*\widehat{\otimes}_{P(n)_*}P(n)^*(X)=
   \inverselimit{}\left(R_*
   \otimes_{P(n)_*}P(n)^*(Z)\right), \] 
where the inverse limit is taken over $\Lambda(X)$.

\begin{lemma}\label{cohomology}
Suppose that $R_*$ is Landweber exact over $P(n)_*$ and even-periodic.
Furthermore, suppose that $R_0$ is Noetherian and linearly compact.
Then the functor $R_*\widehat{\otimes}_{P(n)_*}P(n)^*(-)$
is a complex oriented multiplicative cohomology theory.
\end{lemma}

\proof
Set $R^*(-)=R_*\widehat{\otimes}_{P(n)_*}P(n)^*(-)$.
It is easy to see that $R^*(-)$ takes coproducts to products.
The exactness of $R^*(-)$ follows from 
Lemma~\ref{exact-inverse-limit}
(cf. \cite[Proposition~2.3.16]{Hovey-Palmieri-Strickland}).
The natural transformation $P(n)^*(-)\to R^*(-)$
gives us an orientation of $R^*(-)$.
\qqq

We suppose that $p$ is odd if $n>0$.
In \cite{Wurgler-2} W\"{u}rgler determined the structure
of the co-operation ring $P(n)_*(P(n))$,
which is given as follows:
\[ P(n)_*(P(n))=P(n)_*[t_1,t_2,\ldots ]\otimes
   \Lambda(a_0,a_1,\ldots,a_{n-1}) ,\]
where $|t_i|=2(p^i-1)$ and $|a_i|=2p^i-1$.
In particular, $P(n)_*(P(n))$ is free over $P(n)_*$.
Let $k_*$ be an even-periodic $P(n)_*$-algebra.
Hence we have a degree $0$ formal group law $F_k$.
Note that a $P(n)_*$-algebra homomorphism $\phi: P(n)_*(P(n))\to k_*$
factors through $\overline{\phi}: P(n)_*(P(n))/(a_0,\ldots,a_{n-1})\cong
P(n)_*[t_1,\ldots]$.
It is known that there is a one-to-one correspondence
between a $P(n)_*$-algebra homomorphism
$\overline{\phi}: P(n)_*[t_1,\ldots]\to k_*$ and 
the pair $(f,G)$,
where $G$ is a $p$-typical formal group law of degree $0$
such that
the $p$-series $[p](X)\equiv 0$ mod $(X^{p^n})$, 
and $f$ is an isomorphism 
from $F_k\to G$.
Note that $G$ is the degree $0$ formal group laws
associated with $(\phi\circ\eta_R)_*F_{P(n)}$,
where $\eta_R$ is the right unit of the Hopf algebroid
$P(n)_*(P(n))$, and $F_{P(n)}$ is the degree $-2$
formal group law associated with the complex oriented 
cohomology theory $P(n)^*(-)$.  

Let $h_*$ and $k_*$ be graded rings which are  
Landweber exact over $P(n)_*$ and even-periodic.
We suppose that $h_0$ and $k_0$ are Noetherian and linearly compact.
Set $h^*(-)=h_*\widehat{\otimes}_{P(n)_*}P(n)^*(-)$ and
$k^*(-)=k_*\widehat{\otimes}_{P(n)_*}P(n)^*(-)$.
By Lemma~\ref{cohomology},
$h^*(-)$ and $k^*(-)$ are generalized cohomology theories.
We denote by $h$ and $k$ the representing ring spectra, respectively.  

\begin{proposition}\label{operation-FGL}
We suppose that $p$ is odd if $n>0$.
If $h_*$ and $k_*$ are 
Landweber exact over $P(n)_*$ and even-periodic,
then the map
$\Xi:\mathbf{Mult}(h,k)\longrightarrow
\mathbf{FGL}((F_k,k_0),(F_h,h_0))$
is a bijection.
Furthermore,
if $h_0$ and $k_0$ are Noetherian and linearly compact,
then $\Xi^c:\mathbf{Mult}^c(h,k)\longrightarrow
\mathbf{FGL}^c((F_k,k_0),(F_h,h_0))$
is also a bijection.
\end{proposition}

\proof
For $(f,\alpha)\in \mathbf{FGL}(F_k,F_h)$,
we construct a multiplicative operation $\theta: h^*(-)\to k^*(-)$.
Since $k_*(-)\cong k_* \otimes P(n)_*(-)$,
we have $k_*(P(n))\cong k_*\otimes P(n)_*(P(n))$,
and hence $k_*(P(n))$ is free over $k_*$.
Then $k^0(P(n))\cong \mbox{\rm Hom}_{k_*}(k_*(P(n)),k_*)
\cong \mbox{\rm Hom}_{P(n)_*}(P(n)_*(P(n)),k_*)$.
It is easy to see that 
a multiplicative operation $P(n)^*(-)\to k^*(-)$
corresponds to a $P(n)_*$-algebra homomorphism
$P(n)_*(P(n))\to k_*$.  
Hence we obtain a $P(n)_*$-algebra homomorphism 
$\phi:P(n)_*(P(n))\to k_*$ such that 
$\phi\circ \eta_R$ corresponds to 
$F_h$.
This gives a multiplicative operation
$\varphi: P(n)^*(-)\to k^*(-)$.
Note that $\varphi$ induces $\phi\circ\eta_R$
on the coefficient rings,
and $F_h$ is the degree $0$ formal group law
associated with  $(\phi\circ\eta_R)_*F_{P(n)}$.
By using the ring homomorphism
$\widetilde{\alpha}:h_*\to k_*$,
we may extend $\varphi$ to 
a multiplicative operation
$\theta: h^*(-)=h^*\widehat{\otimes}P(n)^*(-)\to k^*(-)$.
Then it is easy to check that 
this construction gives the inverse of $\Xi$.
If $\alpha: h_0\to k_0$ is continuous,
then $\theta\in\mathbf{Mult}^c(h,k)$ 
by Remark~\ref{cotinuity-remark}.
\qqq


\section{Multiplicative natural transformation $\Theta$}
\label{generalized-Chern-ch}

In this section we suppose that $p$ is an odd prime.
We construct a multiplicative natural transformation 
$\Theta$ from $E^*(-)$ to $K^*(-)\widehat{\otimes}_{\mathbf F}L$,
which is equivariant under the action of $G_{n+1}$.
It is shown that $\Theta$ induces an isomorphism
of ${\mathcal G}$-modules
between $E^*(X)\widehat{\otimes}_{E}L$ and 
$K^*(X)\widehat{\otimes}_{\mathbf F}L$.
This implies that the $G_n$-module $K^*(X)$ is
naturally isomorphic to $H^0(S_{n+1};E^*(X)\widehat{\otimes}_{E}L)$
for all spectra $X$.
Hence we can recover the $G_n$-module structure of $K^*(X)$
from the $G_{n+1}$-module structure of $E^*(X)$.

%

Recall that  ${\mathbf F}$ is an algebraic extension of $\fp{}$
which contains the finite fields $\fp{n}$ and $\fp{n+1}$.
Set $E_*={\mathbf F}\power{u_n}[u^{\pm 1}]$, 
where the degree of $u_n$ is $0$ and the degree of $u$ is $-2$.
Abbreviate to $E$ the degree $0$ subring $E_0={\mathbf F}\power{u_n}$.
We consider that $E_*$ is a $P(n)_*$-algebra
by the ring homomorphism $P(n)_*\to E_*$
given by $v_n\mapsto u_n u^{-(p^n-1)}, 
v_{n+1}\mapsto u^{-(p^{n+1}-1)}, v_i\mapsto 0\ (i>n+1)$.
Then $E_*$ is an even-periodic Landweber exact $P(n)_*$-algebra,
and $E$ is complete Noetherian local ring.  
Hence, by Lemma~\ref{cohomology},
the functor $E^*(-)=E_*\widehat{\otimes}_{P(n)_*}P(n)^*(-)$
is a generalized cohomology theory.
We denote by ${\mathbb E}$ the representing ring spectrum.
Then the degree $0$ formal group law associated with
$E^*(-)$ is $F_{n+1}$.
By Proposition~\ref{operation-FGL},
there is a one-to-one correspondence between
$\mathbf{Mult}^c(\mathbb{E},\mathbb{E})$ and $\mathbf{Aut}^c(F_{n+1},E)$.
By 
Lemma~\ref{automorphism_group_F_n+1},
the automorphism group of $F_{n+1}$ over $E$
is isomorphic to $G_{n+1}$. 
In particular, $G_{n+1}$ acts on the cohomology theory
$E^*(-)$ as multiplicative stable operations.

Let $K_*={\mathbf F}[w^{\pm 1}]$, where $|w|=-2$.
There is a ring homomorphism $P(n)_*\to K_*$
given by $v_n\mapsto w^{-(p^n-1)}, v_i\mapsto 0\ (i>n)$.
By Lemma~\ref{cohomology},
$K^*(-)=K_*\widehat{\otimes}_{P(n)_*}P(n)^*(-)$
is a complex oriented cohomology theory.
We denote by ${\mathbb K}$ the representing ring spectrum. 
Then the associated degree $0$ formal group law
is $H_n$.
By definition, the automorphism group of $(H_n,{\mathbf F})$
is $G_n$.
Hence 
the automorphism group of $K^*(-)$
as multiplicative cohomology theory is 
$G_n$ by Proposition~\ref{operation-FGL}.



The following is the main theorems of this note.

\begin{theorem}\label{main-theorem}
There is a multiplicative stable cohomology operation
\[ \Theta :E^*(-)\to K^*(-)\widehat{\otimes}_{\mathbf F}L \]
such that $\Theta$ is equivariant with respect to
the action of $G_{n+1}$.
Furthermore, $\Theta$ induces an isomorphism
\[ E^*(X)\widehat{\otimes}_{E} L\cong 
   K^*(X)\widehat{\otimes}_{\mathbf F}L,\]
as ${\mathcal G}$-modules for all $X$.
\end{theorem}

\proof
The even-periodic cohomology theory $K^*(-)\widehat{\otimes}_{\mathbf F}L$
is obtained by the even periodic Landweber exact $P(n)_*$-algebra 
$L[w^{\pm 1}]$ given by $v_n\mapsto w^{-(p^n-1)}, v_i\ (i>n)\mapsto 0$.
The associated degree $0$ formal group law is the
Honda group law $H_n$ of height $n$ over $L$.
By Proposition~\ref{operation-FGL},
the automorphism group of $K^*(-)\widehat{\otimes}_{\mathbf F}L$
as a multiplicative cohomology theory is 
$\mathbf{Aut}(H_n, L)$. 
By Theorem~\ref{FGL-main-theorem-v3},
${\mathcal G}$ acts on $K^*(-)\widehat{\otimes}_{\mathbf F}L$
as multiplicative cohomology operations.

By Proposition~\ref{operation-FGL},
$\mathbf{Mult}({\mathbb E},{\mathbb K}\widehat{\otimes}_{\mathbf F}L)\cong
\mathbf{FGL}((H_n,L), (F_{n+1},E))$.
We have the ring homomorphism 
$\alpha: {\mathbf F}\power{u_n}=E\hookrightarrow M\hookrightarrow
L=({\mathbf K}\widehat{\otimes}_{\mathbf F}L)_0$,
and the isomorphism
$\Phi^{-1}: H_n\to \alpha_*F_{n+1}$ over $L$.
Then $(\Phi^{-1},\alpha)\in \mathbf{FGL}((H_n,L),(F_{n+1},E))$
defines a multiplicative natural transformation
$\Theta: E^*(-)\to K^*(-)\widehat{\otimes}_{\mathbf F}L$.
Then $\Theta$ extends to the natural transformation
$\Theta\widehat{\otimes}L: E^*(-)\widehat{\otimes}_{E}L
\to K^*(-)\widehat{\otimes}_{\mathbf F}L$.
Note that $E^*(-)\widehat{\otimes}_{E}L$
is an even-periodic Landweber exact cohomology theory over $P(n)_*$.
Since $\Theta\widehat{\otimes}L$ induces an
isomorphism on the coefficient rings,
$\Theta\widehat{\otimes}L$ is an isomorphism 
of cohomology theories.
Furthermore, by Theorem~\ref{FGL-main-theorem-v3} 
and Proposition~\ref{operation-FGL},
$\Theta\widehat{\otimes}L: 
E^*(X)\widehat{\otimes}_{E}L\stackrel{\cong}{\to}
K^*(X)\widehat{\otimes}_{\mathbf F}L$
is an isomorphism of ${\mathcal G}$-modules,
and $\Theta$ is equivariant under the action of $G_{n+1}$.
%
\qqq

\begin{lemma}
The invariant ring of $L[u^{\pm 1}]$ under the action of
$S_{n+1}$ is $K_*$:
\[ H^0(S_{n+1};L[u^{\pm 1}]) 
   =K_* .\]
\end{lemma}

\proof
Let $M^1_n=v_{n+1}^{-1}BP_*/(p,v_1,\ldots,v_{n-1},v_n^{\infty})$.
By \cite[Theorem~5.10]{MRW},
$\mbox{\rm Ext}^0_{BP_*(BP)}(BP_*,M^1_n)$ is
the direct sum of the finite torsion submodules
and the $K(n)_*/k(n)_*$ generated by $1/v_n^j,\ j\ge 1$ 
as a $k(n)_*$-module.
Then as in \cite[\S5.3]{Torii}
$H^0(S_{n+1};R_*)={\mathbf F}[v_n]$,
where $v_n=u_n u^{-(p^n-1)}$.
By \cite[Lemma~5.9]{Torii},
$H^0(S_{n+1};M_*)$ is the localization of
$H^0(S_{n+1};R_*)$ by inverting the invariant element $v_n$.
Hence $H^0(S_{n+1};M_*)={\mathbf F}[v_n^{\pm 1}]$.

By \cite[Lemma~3.7]{Torii}, 
$w=\Phi_0^{-1}u\in L$
is invariant under the action of $S_{n+1}$.
Let $a$ be a degree $2n$
invariant element in $L[u^{\pm 1}]$.
Then $b=a w^n$
is also invariant.
Let $\phi(X)\in M[X]$ be the minimal polynomial of $b$.
Then $\phi(b)=0$.
Since $b$ is invariant under the action of $S_{n+1}$,
$\phi^g(b)=0$ for all $g\in S_{n+1}$.
Hence $\phi^g(X)$ is also the minimal polynomial of $b$.
This implies that $\phi(X)$ is a polynomial over
$H^0(S_{n+1};M)={\mathbf F}$.
Hence $b\in \overline{{\mathbf F}}\cap L={\mathbf F}$.
This completes the proof. 
\qqq

\begin{corollary}\label{main-corollary}
There are natural isomorphisms of $G_n$-modules\mbox{:}
\[ \begin{array}{rcl}
    K^*(X)&\cong& H^0(S_{n+1};E^*(X)\widehat{\otimes}_{E}
    L),\\[2mm]
    K_*(X)&\cong& H^0(S_{n+1};E_*(X){\otimes}_{E}
    L),\\
   \end{array} \]
for all spectra $X$.
If $X$ is a space, then 
these are also isomorphisms of graded commutative rings.
\end{corollary}

\proof
We have the natural isomorphism of ${\mathcal G}$-modules:
$K^*(X)\widehat{\otimes}_{\mathbf F}L\cong 
E^*(X)\widehat{\otimes}_{E}L$.
The action of the subgroup $S_{n+1}\subset {\mathcal G}$
on the left hand side is obtained from the action on $L$ only.
Hence $H^0(S_{n+1};K^*(X)\widehat{\otimes}_{\mathbf F}L)=K^*(X)$.
This completes the proof of the cohomology case.
The homology case is obtained by the similar way.
\qqq

\section{Lift to characteristic $0$}\label{lift-to-ch-0}

In this section we lift $\Theta$ to
the multiplicative natural transformation $ch$
of the characteristic $0$ cohomology theories.
Then we prove that $ch$ induces a natural isomorphism 
of cohomology theories with stable cohomology operations
if the coefficients are sufficiently extended.
Note that in this section we do not assume that $p$ is an odd prime.

We recall that ${\mathbf F}$ is an algebraic extension
of $\fp{}$ which contains the finite fields $\fp{n}$ and $\fp{n+1}$.
We define graded rings $E_{n,*}$ and $E_{{n+1},*}$ as follows:
\[ \begin{array}{rcccl}
    E_{n,*}&=& E_n[w^{\pm 1}]&
             =&W\power{w_1,\ldots,w_{n-1}}[w^{\pm 1}],\\[2mm]
    E_{{n+1},*}&=&E_{n+1}[u^{\pm 1}]&
              =& W\power{u_1,\ldots,u_n}[u^{\pm 1}],\\
   \end{array}\]
where $W=W({\mathbf F})$ is the ring of Witt vectors with
coefficients in ${\mathbf F}$.
The grading of $E_{n,*}$ is given by
$|w_i|=0\ (1\le i<n)$ and $|w|=-2$,
and the grading of $E_{{n+1},*}$ is given by
$|u_i|=0\ (1\le i\le n)$ and $|u|=-2$.  
Let $r_n:BP_*\to E_{n,*}$ be the ring homomorphism
given by $r_n(v_i)=w_iw^{-(p^i-1)}\ (1\le i<n), 
r_n(v_n)=w^{-(p^n-1)}, r_n(v_i)=0\ (i>n)$,
and let $r_{n+1}: BP_*\to E_{{n+1},*}$ be the ring homomorphism
given by
$r_{n+1}(v_i)=u_i u^{-(p^i-1)}\ (1\le i\le n),
 r_{n+1}(v_{n+1})=u^{-(p^{n+1}-1)}, 
 r_{n+1}(v_i)=0\ (i>n+1)$.
These gives $E_{n,*}$ and $E_{{n+1},*}$
even-periodic Landweber exact $BP_*$-algebra structures.
Hence, by Lemma~\ref{cohomology},
$E_n^*(-)=E_{n,*}\widehat{\otimes}_{BP_*}BP^*(-)$ and 
$E_{n+1}^*(-)=E_{{n+1},*}\widehat{\otimes}_{BP_*}BP^*(-)$
are generalized cohomology theories.
Then there are associated degree $0$ formal group laws 
$\widetilde{F}_n$ and $\widetilde{F}_{n+1}$ 
over $E_n$ and $E_{n+1}$, respectively.
By Lemma~\ref{universal-deformation},
$(\widetilde{F}_n, E_n)$ and
$(\widetilde{F}_{n+1},E_{n+1})$ 
are universal deformations of $(H_n,{\mathbf F})$ and
$(H_{n+1}, {\mathbf F})$, respectively. 

 

Let $R=W\power{u_n}$.
We denote by $S$ the $p$-adic completion of 
$R[u_n^{-1}]$: $S=(W((u_n)))^{\wedge}_p$.
Then $S$ is a complete discrete valuation ring
with uniformizer $p$   
and residue field $M={\mathbf F}((u_n))$. 
In particular,
$S$ is a Henselian ring. 
We recall the following lemma on Henselian rings.

\begin{lemma}[{cf.~\cite[Proposition~I.4.4.]{Milne}}]\label{henselian}
Let $A$ be a Henselian ring with residue field $k$.
Then the functor $B\mapsto B\otimes_A k$ induces
an equivalence between the category of finite \'{e}tale  
$A$-algebras and the category of finite \'{e}tale $k$-algebras.
\end{lemma}

In \cite{Torii},
we have constructed a sequence of finite separable extensions of $M$:
\[ M=L_{-1}\to L_0\to L_1\to\cdots,\]
where $L_i$ is obtained by adjoining the coefficients 
$\Phi_0,\Phi_1,\ldots,\Phi_i$ of the isomorphism 
$\Phi: F_{n+1}\stackrel{\cong}{\to} H_n$. 
By definition,
$L=\directlimit{i}L_i=\cup_i L_i$ and
we have shown that $L_i$ is stable under the action of ${\mathcal G}$
for all $i$.
By Lemma~\ref{henselian},
we obtain a sequence of finite \'{e}tale 
$S$-algebras:
\[ S=S_{-1}\to S_0\to S_1\to \cdots .\] 
We denote by $S_{\infty}$ the direct limit $\directlimit{i}S_i$
and $T$ the $p$-adic completion of $S_{\infty}$.

\begin{lemma}       
The ring $T$ is a complete discrete valuation ring 
of characteristic $0$ 
with uniformizer $p$ and residue field $L=\directlimit{i}L_i$.
\end{lemma}

\proof
Since $L_i$ is a separable extension over $M$,
we can take $a\in L_i$ such that $L_i=M(a)$.
Let $f(X)\in M[X]$ be the minimal polynomial of $a$
and $\widetilde{f}(X)\in S[X]$ a monic polynomial
which is a lift of $f(X)$.
Then $S_i\cong S[X]/(\widetilde{f}(X))$.
Then we see that $S_i$ is a complete 
discrete valuation ring with uniformizer $p$
and residue field $L_i$.
This implies that 
$S_{\infty}$ is also a discrete valuation ring with uniformizer $p$
and residue field $L$.
Then the lemma follows from the fact that
$T$ is the $p$-adic completion of $S_{\infty}$.
\qqq

We abbreviate $T\power{w_1,\ldots,w_{n-1}}$
and $T\power{u_1,\ldots,u_{n-1}}$ by
$T\power{w_i}$ and $T\power{u_i}$, respectively, etc.
Then we obtain a sequence of finite \'{e}tale $S\power{u_i}$-algebras:
\[ S\power{u_i}= S_{-1}\power{u_i}
               \to S_0\power{u_i}
               \to S_1\power{u_i}\cdots, \]
and 
$T\power{u_i}$ is the $I_n$-adic
completion of 
$\directlimit{j}S_j\power{u_i}$,
where $I_n=(p,u_1,\ldots,u_{n-1})$.

The ring homomorphisms
$BP_*\to E_{n,*}\hookrightarrow 
T\power{w_i}[w^{\pm 1}]$
and
$BP_*\to E_{{n+1},*}\hookrightarrow 
T\power{u_i}[u^{\pm 1}]$
satisfy the Landweber exact condition.
Also $T\power{w_i}[w^{\pm 1}]$ and $T\power{u_i}[u^{\pm 1}]$
are even-periodic, and 
the degree $0$ subring $T\power{w_i}$ and $T\power{u_i}$
are complete Noetherian local rings.
By Lemma~\ref{cohomology},
the following two functors are generalized cohomology theories:
\[ \begin{array}{rcl}
   E_n^*(X)\widehat{\otimes}_{E_n}T\power{w_i}&:=&
   \inverselimit{\Lambda(X)}(E_n^*(X_{\alpha})\otimes_{E_n}T\power{w_i}),\\
   E_{n+1}^*(X)\widehat{\otimes}_{E_{n+1}}T\power{u_i}&:=&
   \inverselimit{\Lambda(X)}
   (E_{n+1}^*(X_{\alpha})\otimes_{E_{n+1}}T\power{u_i}).\\
   \end{array}\]
The degree $0$ formal group laws associated with
$E_n^*(-)\widehat{\otimes}T\power{w_i}$ and 
$E_{n+1}^*(-)\widehat{\otimes}T\power{u_i}$
are $(\widetilde{F}_n,T\power{w_i})$ and 
$(\widetilde{F}_{n+1},T\power{u_i})$, respectively.

\begin{lemma}\label{universality}
The formal group laws $(\widetilde{F}_{n+1},T\power{u_i})$
and $(\widetilde{F}_n,T\power{w_i})$ are
universal deformations of
$(F_{n+1},L)$ and $(H_n,L)$,respectively,
on the category of complete Noetherian local $T$-algebras.
\end{lemma}

\proof
From the fact that $(\widetilde{F}_n,E_n)$ is a universal
deformation of $(H_n,{\mathbf F})$,
it is easy to see that 
$(\widetilde{F}_n, T\power{w_i})$ is a universal deformation of
$(H_n,L)$.
From the form of the $p$-series of $\widetilde{F}_{n+1}$ given by 
(\ref{p-series-tilde-Fn+1}),
we see that $(\widetilde{F}_{n+1},T\power{u_i})$
is a universal deformation of $(F_{n+1},L)$.
\qqq


\begin{corollary}
The action of $G_{n+1}$ on $(\widetilde{F}_{n+1},E_{n+1})$
extends to an action on $(\widetilde{F}_{n+1},T\power{u_i})$
such that the induced action on $(F_{n+1},L)$ 
coincides with the action of Theorem~\ref{FGL-main-theorem-v3}.
\end{corollary}

\proof
It is sufficient to show that 
the action of $G_{n+1}$ on $E_{n+1}$ extends
to an action on $T\power{u_i}$.
For $g\in G_{n+1}$,
$u_n^g$ is a unit multiple of $u_n$ 
modulo $(p,u_1,\ldots,u_{n-1},u_n^2)$.
Hence the ring homomorphism
$E_{n+1}\stackrel{g}{\to}E_{n+1}\to
 (E_{n+1}[u_n^{-1}])^{\wedge}_{I_n}=S\power{u_i}$
extends to a ring homomorphism
$E_{n+1}[u_n^{-1}]\to S\power{u_i}$
This induces a ring homomorphism 
$S\power{u_i}\to S\power{u_i}$
and defines an action of $G_{n+1}$ on $S\power{u_i}$.
Since $S_j\power{u_i}\to S_{j+1}\power{u_i}$
is \'{e}tale for $j\ge -1$
and $L_j$ is stable under the action of $G_{n+1}$ on $L$,
the action on $S\power{u_i}$ extends to 
$S_j\power{u_i}$ uniquely and compatibly by Lemma~\ref{henselian}.
Hence we obtain an action on $\directlimit{j}(S_j\power{u_i})$
and its $I_n$-adic completion $T\power{u_i}$.
\qqq


We denote the action of $G_{n+1}$ on $(\widetilde{F}_{n+1},T\power{u_i})$ 
by $\Upsilon(g)=(t(g), \upsilon(g)): 
(\widetilde{F}_{n+1},T\power{u_i})\to 
(\widetilde{F}_{n+1},T\power{u_i})$ for $g\in G_{n+1}$.

\begin{corollary}
The $(n+1)$th extended Morava stabilizer group 
$G_{n+1}$ acts on the cohomology theory
$E_{n+1}^*(-)\widehat{\otimes}T\power{u_i}$ 
as multiplicative cohomology operations.
\end{corollary}

\proof
This follows from Proposition~\ref{operation-FGL}. 
\qqq

Recall that $S_n$ and $G_n$ are 
identified with the Galois groups $\gal{L/M}$ and
$\gal{L/\fp{}((u_n))}$, respectively,
through the action of $G_n$ on $L$
(\cite{Gross, Torii}). 

\begin{lemma}
The action of $G_n$ on $L$ lifts to the action on $T$.
\end{lemma}

\proof
Since $L_i$ is stable under the action of $G_n$ on $L$
for all $i\ge -1$,
the action of $G_n$ on $L_i$ lifts to the action on
$S_i$ compatibly by Lemma~\ref{henselian}.
This induces an action on $S_{\infty}$.
Since $T$ is the $p$-adic completion of $S_{\infty}$,
we obtain an action on $T$
which is a lift of the action on $L$.
\qqq

We denote this action of $G_n$ on $T$ by $\tau(g): T\to T$
for $g\in G_n$.  
Since the actions of $G_n$ on $E_n$ and $T$
are compatible on $W$,
the diagonal action defines an action of $G_n$ on $T\power{w_i}=
T\widehat{\otimes}_W W\power{w_i}$. 
Then we obtain an extension of the action of $G_n$ 
on $(\widetilde{F}_n, E_n)$ to $(\widetilde{F}_n, T\power{w_i})$.
We denote this action of $G_n$ on $(\widetilde{F}_n, T\power{w_i})$
by $\Omega(g)=(s(g),\omega(g)):
(\widetilde{F}_n, T\power{w_i})\to (\widetilde{F}_n, T\power{w_i})$
for $g\in G_n$.

\begin{corollary}
The $n$th extended Morava stabilizer group
$G_n$ acts on 
$E_n^*(-)\widehat{\otimes} T\power{w_i}$ as 
multiplicative cohomology operations. 
\end{corollary}

\proof
This follows from Proposition~\ref{operation-FGL}. 
\qqq

\begin{lemma}
There is a unique isomorphism 
$(\widetilde{\Phi},\widetilde{\varphi}):
(\widetilde{F}_{n+1},T\power{u_i})\to
(\widetilde{F}_n,T\power{w_i})$
such that $\widetilde{\varphi}$ is a continuous $T$-algebra
homomorphism and $\widetilde{\Phi}$ induces
$\Phi$ on the residue fields.
\end{lemma}

\proof
Since there is an isomorphism
$(\Phi,id_L):(\widetilde{F}_{n+1},L)\to (H_n, L)$, 
the lemma follows from Lemma~\ref{key-lemma-ch0}.
\qqq

%

\begin{lemma}\label{fundamental-commutativity}
For $g\in G_n$
there is a commutative diagram:
\[ \begin{array}{ccc}
    (\widetilde{F}_{n+1},T\power{u_i}) & 
    \stackrel{(X,\theta(g))}{\longrightarrow} &
    (\widetilde{F}_{n+1},T\power{u_i}) \\
    {\scriptstyle (\widetilde{\Phi},\widetilde{\varphi})}
    \bigg\downarrow & & \bigg\downarrow 
    {\scriptstyle (\widetilde{\Phi},\widetilde{\varphi})}\\
    (\widetilde{F}_n,T\power{w_i}) & 
    \stackrel{\Omega(g)}{\longrightarrow} &
    (\widetilde{F}_n,T\power{w_i}), \\
   \end{array}\]
where $\theta(g):T\power{u_i}\to T\power{u_i}$
is given by $\theta(g)(t)=\tau(g)(t)$ for $t\in T$
and $\theta(g)(u_i)=u_i$ for $1\le u_i<n$.
\end{lemma}

\proof
Note that $(\theta(g)\circ\widetilde{\varphi})|_T=
\tau(g)=(\widetilde{\varphi}\circ \omega(g))|_T$.
The diagram induced on the residue field is commutative
by definition of the action of $G_n$ on $(F_{n+1},L)\cong (H_n,L)$.
Then the lemma follows from the universality of
$(\widetilde{F}_n,T\power{w_i})$.
\qqq

\begin{corollary}\label{commutativity-actions}
The pro-finite group ${\mathcal G}$ acts on 
$(\widetilde{F}_{n+1},T\power{u_i})
\cong (\widetilde{F}_n,T\power{w_i})$ such that
the action of the subgroup $G_{n+1}$ coincides with
$\Upsilon$,
and the action of the subgroup $G_n$ coincides with 
$\Omega$.
\end{corollary}

\proof
We have the action $\Upsilon$ of $G_{n+1}$ on 
$(\widetilde{F}_{n+1},T\power{u_i})$ and 
the action $\Omega$ of $G_n$ on $(\widetilde{F}_n,T\power{w_i})$.
The action of the subgroup $\Gamma$ of $G_{n+1}$ 
on $(\widetilde{F}_{n+1},T\power{u_i})$ coincides with
the action on $(\widetilde{F}_n,T\power{w_i})$
as the subgroup of $G_n$ under the isomorphism
$(\widetilde{\Phi},\widetilde{\varphi})$.
Hence it is sufficient to show that 
the following diagram commutes for
$g\in S_{n+1}$ and $h\in S_n$:
\[  \begin{array}{ccc}
     (\widetilde{F}_{n+1},T\power{u_i}) &
     \stackrel{(X,\theta(h))}{\longrightarrow} &
     (\widetilde{F}_{n+1},T\power{u_i}) \\
     {\scriptstyle \Upsilon(g)}
     \bigg\downarrow & & \bigg\downarrow 
     {\scriptstyle \Upsilon(g)} \\
     (\widetilde{F}_{n+1},T\power{u_i}) &
     \stackrel{(X,\theta(h))}{\longrightarrow} &
     (\widetilde{F}_{n+1},T\power{u_i}). \\ 
    \end{array}\]
Note that the induced diagram on the residue field $L$
commutes.

Since $u_n^g\in E_{n+1}\subset T\power{u_i}$,
$(\theta(h)\circ \upsilon(g))(u_n)=u_n^g=
(\upsilon(g)\circ \theta(h))(u_n)$.
Hence $(\theta(h)\circ \upsilon(g))|_S=
(\upsilon(g)\circ \theta(h))|_S$.
From the fact that $S_i$ is an \'{e}tale $S$-algebra,
$T\power{u_i}$ is complete, and the induced homomorphisms
on the residue field coincide,
we see that $(\theta(h)\circ \upsilon(g))|_{S_i}=
(\upsilon(g)\circ \theta(h))|_{S_i}$ for all $i$.
Hence $(\theta(h)\circ \upsilon(g))|_{S_{\infty}}=
(\upsilon(g)\circ \theta(h))|_{S_{\infty}}$ and 
$(\theta(h)\circ \upsilon(g))|_T=
(\upsilon(g)\circ \theta(h))|_T$.
Then the corollary follows from the universality of
$(\widetilde{F}_{n+1},T\power{u_i})$.
\qqq

\begin{theorem}\label{Characteristic-zero}
There is a multiplicative stable cohomology operation
\[ ch :E_{n+1}^*(-)\to E_n^*(-)\widehat{\otimes}_{E_n} T\power{w_i} \]
such that $ch$ is equivariant with respect to
the action of $G_{n+1}$.
Furthermore, $ch$ induces a natural isomorphism
\[ E_{n+1}^*(X)\widehat{\otimes}_{E_{n+1}} T\power{u_i}\cong 
   E_n^*(X)\widehat{\otimes}_{E_n} T\power{w_i},\]
as ${\mathcal G}$-modules for all spectra $X$.
\end{theorem}

\proof
As in the proof of Theorem~\ref{main-theorem},
this follows from Lemma~\ref{commutativity-actions}
and Proposition~\ref{operation-FGL}.
\qqq

\begin{remark}\rm
As in Lemma~\ref{fundamental-commutativity},
we can show that the following diagram commutes for $g\in G_{n+1}$:
\[ \begin{array}{ccc}
    (\widetilde{F}_{n+1}, T\power{u_i}) & 
    \stackrel{\Upsilon(g)}{\longrightarrow} & 
    (\widetilde{F}_{n+1}, T\power{u_i}) \\
    {\scriptstyle (\widetilde{\Phi},\widetilde{\varphi})}
    \bigg\downarrow & & \bigg\downarrow 
    {\scriptstyle (\widetilde{\Phi},\widetilde{\varphi})}\\
    (\widetilde{F}_n, T\power{w_i}) & 
    \stackrel{(X,\mu(g))}{\longrightarrow} &
    (\widetilde{F}_n, T\power{w_i}),\\
   \end{array}\] 
where $\mu(g)$ is given by 
$\mu(g)(t)=\widetilde{\varphi}^{-1}(\upsilon(g)(t))$ for $t\in T$
and $\mu(g)(w_i)=w_i$ for $1\le i<n$. 
This implies that 
there is a $G_n$-equivariant natural homomorphism
\[ E_n^*(X)\longrightarrow 
   H^0(S_{n+1};E_{n+1}^*(X)\widehat{\otimes} T\power{u_i}) ,\]
which is a homomorphism of graded commutative rings if
$X$ is a space.
\end{remark}



\begin{flushleft}
Department of Mathematics, Okayama University, \\
Okayama 700--8530, Japan\\
{\it E-mail}\,: torii@math.okayama-u.ac.jp\\
\end{flushleft}

{\footnotesize

}

\end{document}